\documentclass{article}
\usepackage{amsthm}
\usepackage{amsfonts,amssymb,amsmath}
\usepackage{enumitem}
\usepackage{titlesec}
\usepackage{fouriernc}
\usepackage[T1]{fontenc}
\usepackage{chngcntr}
\counterwithin{figure}{section}

\newlist{gcases}{enumerate}{1}
\setlist[gcases,1]{
  label={{\it Case}~{\it \Alph*}.},
  topsep=0ex,
  leftmargin=0in,
  labelsep=.1in,
  itemindent=.7in,
  itemsep=0ex 
}
\newlist{tenumerate}{enumerate}{1}
\setlist[tenumerate,1]{
  label={(\arabic*)},
  topsep=0ex,
  leftmargin=.3in,
  labelsep=.1in,
  itemindent=0in,
  itemsep=0ex
}

\newlist{titemize}{enumerate}{1}
\setlist[titemize,1]{
  label={$\bullet$},
  topsep=0ex,
  leftmargin=.3in,
  labelsep=.1in,
  itemindent=0in,
  itemsep=0ex
}

\titleformat{\section}
  {\normalfont\bfseries}   
  {}                             
  {0pt}                          
  {Section \thesection\quad}    

\titleformat{name=\section,numberless}
  {\normalfont\bfseries}
  {}
  {0pt}
  {}

\newlength{\tabwidth}
\newlength{\tabheight}
\newlength{\tabrule}
\newlength{\tabwidthx}
\newlength{\tabheightx}

\def\gentabbox#1#2#3#4{\vbox to \tabheight{\setlength{\tabrule}{#3}%
  \setlength{\tabwidthx}{#1\tabwidth}\addtolength{\tabwidthx}{\tabrule}%

\setlength{\tabheightx}{#2\tabheight}\addtolength{\tabheightx}{-\tabheight}%
  \hbox to #1\tabwidth{%
 \hspace{-0.5\tabrule}\rule{\tabrule}{#2\tabheight}\hspace{-\tabrule}%
    \vbox to #2\tabheight{\hsize=\tabwidthx%
      \vspace{-0.5\tabrule}\hrule width\tabwidthx height\tabrule%
      \vspace{-0.5\tabrule}\vfil%
      \hbox to \tabwidthx{\hss#4\hss}%
        \vfil\vspace{-0.5\tabrule}%
      \hrule width\tabwidthx height\tabrule\vspace{-0.5\tabrule}}%
 \hspace{-\tabrule}\rule{\tabrule}{#2\tabheight}\hspace{-0.5\tabrule}}%
  \vspace{-\tabheightx}}}
\def\genblankbox#1#2{\vbox to \tabheight{\vfil\hbox to
#1\tabwidth{\hfil}}}
\def\tabbox#1#2#3{\gentabbox{#1}{#2}{0.4pt}{\strut #3}}

\catcode`\:=13 \catcode`\.=13 \catcode`\;=13
\catcode`\>=13 \catcode`\^=13
\def:#1\\{\hbox{$#1$}}
\def.#1{\tabbox{1}{1}{$#1$}}
\def>#1{\tabbox{2}{1}{$#1$}}
\def^#1{\tabbox{1}{2}{$#1$}}
\def;{\genblankbox{1}{1}\relax}
\catcode`\:=12 \catcode`\.=12 \catcode`\;=12
\catcode`\>=12 \catcode`\^=7

\newcommand\T{\mathbf{T}}

\begin{document}

\noindent
{\bf \large Orbital varieties in types $B$ and $C$}

\vspace{.3in}
\noindent
WILLIAM~M.~MCGOVERN \\
\vspace{.1in}
{\it \small Department of Mathematics, Box 354350, University of Washington, Seattle, WA, 98195}

\section*{Introduction}
The purpose of this paper is to correct the proof of the main result of \cite{M99}, which parametrizes orbital varieties in the classical cases in terms of domino tableaux and gives an algorithm for computing the tableau of the orbital variety corresponding to a classical Weyl group element.  We treat types $B$ and $C$ in this paper; a later one will deal with type $D$.  As in \cite{M99}, we will rely heavily on three papers of Garfinkle \cite{G90,G92,G93} and one of the author and Pietraho \cite{MP16}, which use standard domino tableaux to classify Kazhdan-Lusztig cells in types $B,C$, and $D$ .

\section*{Main result}
We recall the main result of \cite{M99}, using the notation and definitions of that paper.

\newtheorem*{thm}{Theorem}
\begin{thm}
Assume that the simple Lie algebra $\mathfrak g$ is of type $X=B$ or $C$.  Orbital varieties contained in the orbit $\mathcal O_{\mathbf p}$ with partition $\mathbf p$ are parametrized by standard domino tableaux of shape $\mathbf p$.  For any $w$ in the Weyl group $W$, the tableau parametrizing the corresponding variety $V(w)$ may be obtained from the left tableau $T_L(w)$ attached by Garfinkle to $w$ by moving through certain open cycles of type $X$ or $D$ whose holes and corners lie in rows of odd length in type $C$, or even length in type $B$, lowering the partition corresponding to the shape of the resulting tableau in the standard dominance order throughout, until this shape becomes an $X$-partition $\mathbf q$.  The orbit $\mathcal O_w$ then coincides with $\mathcal O_{\mathbf q}$.  
\end{thm}

\begin{proof}
As in \cite{M99}, we write $w_1\sim w_2$ whenever $w_,w_2\in W$ lie in the same geometric left cell (so that $V(w_1) = V(w_2))$.  Let $S$ be a subset of the set $\Delta$ of simple roots.  Writing any $w\in W$ uniquely as $w_Sw'$ for some $w_S\in W_S$, the parabolic subgroup of $W$ generated by the reflections in $S$, and $w'$ a minimal length coset representative of $W_S$ in $W$, we argue as in \cite{M99} to show that $v\sim w$ relative to $W$ whenever $v'=w'$ and $v_S\sim w_S$ relative to $W_S$.  

In particular whenever $\alpha,\beta$ are adjacent simple roots of the same length, the operator $T_{\alpha\beta}^R$ defined in \cite{G92} preserves geometric left cells wherever it is defined.  If $\alpha,\beta$ are adjacent but have different lengths, then we need a single-valued truncation $V_{\alpha\beta}^R$ of this operator.  Write $(\T_1',\T_2')$ for the image of $V_{\alpha\beta}^R$ applied to a pair $(\T_1,\T_2)$ on which it is defined.  In type $C$, $V_{\alpha\beta}^R$ is defined on $(\T_1,\T_2)$ if the 1- and 2-dominos of $\T_2$ form a subtableau of shape either $(3,1)$ or $(2,2)$; in the former case we obtain $(\T_1',\T_2')$  from $(\T_1,\T_2)$ by moving both tableaux through the extended cycle of the 2-domino in $\T_2$ relative to $\T_1$ and then transposing the 1- and 2-dominos in the $2\times 2$ box they now occupy in the right tableau.  In the latter case we set $\T_1'=\T_1$ and obtain $\T_2'$  from $\T_2$ by transposing the 1- and 2-dominos in the $2\times 2$ box they occupy (while leaving all other dominos unchanged).  In type $B, V_{\alpha\beta}^R$ is defined if the 0 square and 1- and 2-dominos in $\T_2$ form a subtableau of shape $(3,2)$ or $(3,1,1)$; in the former case we move $(\T_1,\T_2)$ through the extended cycle of the 2-domino in $\T_2$ relative to $\T_1$ and then interchange the positions of the 1- and 2-dominos in the right tableau, while in the latter we just interchange the positions of the 1- and 2-dominos in the right tableau.  Using \cite[2.1.4,2.3.4]{G92} and the above paragraph, we see that $V_{\alpha\beta}^R$ preserves the geometric left cell whenever it is defined.  

We also need an analogous operator $V_D^R$, which is again a single-valued truncation of a larger operator, this time the operator $T_{\alpha_1',C}^R$ defined in \cite[4.4.10]{MP16}.  In type $C$, this operator is defined on pairs $(\T_1,\T_2)$ for which the first four dominos in $\T_2$ form a subtableau with shape $(4,3,1)$ and the 2-domino is vertical (thus in the first column).  Its image is then obtained by moving $(\T_1,\T_2)$ through the extended cycle {\sl in type $D$ } of the 4-domino in $\T_2$ relative to $\T_1$ and then transposing the positions of the 2- and 4-dominos in the $2\times 2$ box they occupy in $\T_2$.  In type $B$, this operator is defined on pairs $(\T_1,\T_2)$ for which the 0 square and first three dominos in $\T_2$ form a subtableau of shape $(4,2,1)$; its image is obtained by moving $(\T_1,\T_2)$ through the extended cycle in type $D$ of the 3-domino in $\T_2$ relative to $\T_1$ and then transposing the positions of the 2- and 3-dominos in the right tableau in the $2\times 2$ box they occupy in $\T_2$. As before, we check that this operator preserves the geometric left cell whenever it is defined.

We now argue as in \cite[p. 2987]{M99}, with a small modification.   In type $C$, given $w\in W$ whose left tableau $T_L(w)$ is not already an $X$-partition, look at the largest odd part $p$ occurring with odd multiplicity in the partition corresponding to the shape of $T_L(w)$.  If the lowest row of length $p$ in this shape is an odd-numbered row, then it is easy to construct a tableau $\T$ having the same shape as $T_R(w)$ with the $m(m+1)$- and $(m+1)^2$-dominos having the same orientation as in $T_R(w)$, for all $m\ge2$ such that dominos with these labels appear in $T_R(w)$, and such that the first two dominos in $\T$ form a subtableau of shape $(3,1)$; moreover, the proof of the main result of \cite{M20} shows that we may apply a sequence $\Sigma$ of operators $T_{\alpha\beta}^R$ and $V_{\alpha\beta}^R$ to $T_R(w)$ in such a way that tableau shapes are preserved throughout and the resulting tableau is $\T$.  Applying the operator $V_{\alpha\beta}^R$ one more time, this time lowering the shape in the dominance order, we arrive at $w'\in W$ in the same geometric left cell as $w$ whose left tableau $T_L(w')$ has lower shape than that of $T_L(w)$.  If instead the lowest row of length $p$ in the shape as $T_R(w)$ is an even row, we argue similarly, producing a tableau $\T$ in the domain of $V_D^R$ whose first four dominos form a subtableau of shape $(4,3,1)$, such that some sequence $\Sigma$ of operators $T_{\alpha\beta}^R$ and $V_{\alpha\beta}^R$ applied to $(T_L(w),T_R(w))$ preserves tableau shapes throughout and transforms the right tableau to $\T$; then applying $V_D^R$ we again get $w'\in W$ in the same geometric left cell as $w$ whose left tableau $T_L(w')$ has a lower shape than that of $T_L(w)$.  Iterating this process, we eventually replace $w$ by an element $v$ in the same geometric left cell the shape of whose left tableau $T_L(v)$ is a $C$-partition.  Further applications of $V_{\alpha\beta}^R$ and $V_D^R$ cannot lower the shape of this partition in the dominance order, since for example if the first two dominos of $T_L(v)$ form a subtableau of shape $(3,1)$, then the cycle of the 2-domino is unboxed in the sense of \cite{G90}; if it were open, then it would be a down cycle and moving through this cycle would produce a tableau whose row lengths do not weakly decrease as one moves down, a contradiction.  In type $B$, we argue similarly, interchanging the roles of even and odd parts, and using the operators $T_{\alpha\beta}^R$ and $V_D^R$ as defined in that case.  

We are therefore reduced to the case where $T_L(w)$ is an $X$-partition.  Now we recall that orbital varieties, like primitive ideals of a fixed regular infinitesimal character, have $\tau$-invariants.  Given a pair $\alpha\beta$ of adjacent simple roots there is a wall-crossing operator $T_{\alpha\beta}^V$ defined on orbital varieties having one of $\alpha,\beta$ in their $\tau$-invariants but not the other \cite[9.9,9.11]{J84}.  If $\alpha$ and $\beta$ have the same length then the image of this operator, say on a variety $V$ with $\alpha$ in its $\tau$-invariant but not $\beta$, is a variety with $\beta$ in its $\tau$-invariant but not $\alpha$, obtained by taking the image of $V$ under the minimal nonsolvable parabolic subgroup $P_\beta$ corresponding to $\beta$ and then intersecting with the nilradical $\mathfrak n_\beta$ of its Lie algebra.  In types $B$ and $C$ we have $T_{\alpha\beta}^V V(w) = V(T_{\alpha\beta}^R(w))$ for $\alpha,\beta$ of the same length \cite[9.11]{J84}, defining $T_{\alpha\beta}^R(w)$ as in \cite{G92}.  If instead $\alpha$ and $\beta$ have different lengths, then the operator $T_{\alpha\beta}^V$ is defined by the same recipe, but behaves in a slightly more complicated way.  It is defined on a variety $V$ having $\alpha$ in its $\tau$-invariant but not $\beta$ {\sl if} some component of $P_\beta(V)\cap\mathfrak n_{\beta}$ has the same dimension as $V$ and can take one or two values in that case (that is, there can be one or two components in the intersection of the same dimension as $V$); these values are then again those of $V(T_{\alpha\beta}^R(w))$.  Specifically, in type $C$,  if the first two dominos of $T_L(w)$ form a subtableau of shape $(2,1,1)$, so that the simple root $\beta=e_2-e_1$ lies in the $\tau$-invariant $\tau(w)$ but the root $\alpha=2e_1$ does not, and if the cycle of the 2-domino in $T_L(w)$ is open, then moving through this cycle raises the shape of $T_L(w)$ in the dominance order and $T_{\alpha\beta}^V$ is not defined on $V(w)$.  Similarly,  if the first two dominos of $T_L(w)$ are vertical and form a $2\times 2$ box, then one of the images of $T_{\beta\alpha}^R(w)$, namely the one with left tableau obtained from $T_L(w)$ by transposing the positions of its first dominos within their $2\times 2$ box and then moving the 2-domino through its cycle, does not appear as an image of $T_{\beta\alpha}(w)$ if this 2-cycle is open.   Otherwise, this map coincides with $T_{\alpha\beta}^R$.  A similar situation holds in type $B$.  We use the maps $T_{\alpha\beta}^V$ to define the notion of generalized $\tau$-invariant of a variety as in \cite{G93} and argue as in that paper to show that the generalized $\tau$-invariant of a variety is a complete invariant, so that two varieties in the same Lie algebra with the same generalized $\tau$-invariant coincide.

Now it is easy to complete the proof.  Given two Weyl group elements $w,w'$ with the same left tableau $T_L(w)=T_L(w')$ having the shape of an $X$-partition, any operator $T_{\alpha\beta}^V$ defined on $V(w)$ is also defined on $V(w')$ and preserves the shapes of $T_L(w)$ and $T_L(w')$, whence in fact it fixes these tableaux, so that $V(w)$ and $V(w')$ have the same generalized $\tau$-invariant and must coincide.  Hence there are at most as many distinct orbital varieties as there are standard domino tableaux with shape an $X$-partition; but the number of such tableaux is the sum of the dimensions of the Springer representations attached to the corresponding nilpotent orbits in the ambient Lie algebra, which is well known to be the number of orbital varieties, so we get a natural bijection between orbital varieties and standard domino tableaux with shape that of an $X$-partition.  Since we know by \cite[pp. 2985,6]{M99} that given an $X$-partition $\mathbf p$ at least one $w\in W$ has left tableau $T_L(w)$ of shape $\mathbf p$ and orbital variety $V(w)$ lying in the orbit $\mathcal O_{\mathbf p}$, the bijection between standard domino tableaux and orbital varieties must send tableaux with shape $\mathbf p$ to varieties lying in $\mathcal O_{\mathbf p}$, as claimed.
\end{proof}

\noindent The above argument also validates the proof of \cite[Theorem 2]{M99}, which asserts that one component of the variety $\mathcal V(w)$ of the simple highest weight module $L(w\gamma)$ of highest weight $w\gamma-\rho$, where $\gamma$ is an antidominant regular integral weight and $\rho$ is half the sum of the positive roots, is the one corresponding to the tableau $S_R(w)$, the unique tableau of special shape obtained from the right tableau $T_R(w)$ by moving through open $X$-cycles (keeping $X$ fixed this time and allowing either up or down cycles).

\end{document}